\documentclass[11pt]{amsart}

\usepackage[margin=1in]{geometry}
\usepackage{amsmath,amssymb,amsthm,mathtools}
\usepackage[colorlinks=true,linkcolor=blue,citecolor=blue,urlcolor=blue]{hyperref}

\newtheorem{theorem}{Theorem}
\newtheorem{lemma}{Lemma}

\newtheorem{remark}{Remark}

\newcommand{\C}{\mathbb C}
\newcommand{\Pb}{\mathbb P}
\newcommand{\R}{\mathbb R}
\newcommand{\Sb}{\mathbb S}

\newcommand{\Z}{\mathbb Z}
\newcommand{\bbone}{\mathbf 1}

\title{A construction of spherical $5$-designs with $O(d^2)$ points}

\author{A. Arman}
 \address{Department of Mathematics, University of Manitoba, Winnipeg, MB, R3T 2N2, Canada}
 \email{andrew0arman@gmail.com}
 
\author{A.\ Bondarenko}
 \address{Department of Mathematical Sciences, Norwegian University of Science and
 	Technology, NO-7491 Trondheim, Norway}
 \email{andriybond@gmail.com}
 \thanks{The second author was supported in part by Grant 334466 of the Research Council of Norway.}

\author{A.\ Prymak}
 \address{Department of Mathematics, University of Manitoba, Winnipeg, MB, R3T 2N2, Canada}
 \email{prymak@gmail.com}
 \thanks{The third author was supported by NSERC of Canada Discovery Grant RGPIN-2026-06488.}

\author{D.\ Radchenko}
 \address{Institut des Hautes \'Etudes Scientifiques, CNRS, Laboratoire Alexandre Grothendieck, 35 route de Chartres, Bures-sur-Yvette 91440, France}
 \email{danradchenko@gmail.com}
\thanks{The fourth author was supported by ERC Starting Grant No. 101078782.}

\date{}


\begin{document}
\maketitle

\begin{abstract}
For every $d\geq1$ we give an explicit equal-weight spherical $5$-design in
$\Sb^{d-1}\subset\R^d$ with at most $72d^2$ points. Our approach utilizes recent construction of complex projective $2$-designs based on Sidon sets.
\end{abstract}

\section{Introduction}

A spherical $t$-design is a finite set of points on a sphere whose average
agrees with the sphere average for every polynomial of degree at most $t$.  Thus
$X\subset \Sb^{d-1}$ is a spherical $t$-design if
\[
        \frac1{|X|}\sum_{x\in X}p(x)=\int_{\Sb^{d-1}}p(x)\,d\sigma_d(x)
\]
for every $d$-variable polynomial $p$ of degree at most $t$, where $\sigma_d$ is the
normalized surface measure.  The concept of spherical designs was introduced by Delsarte--Goethals--Seidel \cite{DGS1977}.

The lower bound of Delsarte--Goethals--Seidel states that a spherical $5$-design in dimension~$d$ must have
 	\[N\ge d(d+1)\]
points. This bound is attained for $d=2,3,7,23$~\cite[Example 8.3]{DGS1977}. Misawa~\cite{Misawa2026} recently gave an explicit construction of spherical $5$-designs in every dimension with $O(d^3)$ points, using lifts from tight fusion frames.

The aim of this note is to give a simple construction of spherical $5$-designs with $O(d^2)$ points in all dimensions. The key ingredient is a method for generating complex projective $2$-designs with $O(d^2)$ points using Sidon sets. The argument we use is essentially the same as the one used by Iosue, Mooney, Ehrenberg and Gorshkov in \cite[Sec. 4.1]{IosueMooneyEhrenbergGorshkov2024}. The earlier use of Sidon sets (more specifically, Singer difference sets) in these problems is due to Bodmann and Haas~\cite{BodmannHaas2016}; they used them to construct weighted complex projective $2$-designs.

\begin{theorem}\label{thm:main}
For every integer $d\geq1$ there is an explicit equal-weight spherical
$5$-design $X_d\subset S^{d-1}\subset\R^d$ with 
\[
        |X_d|\leq 72d^2;
\]
moreover,
\[
|X_d|\le \begin{cases}
    12d(d+2), &\text{ if $d\ge4$ is even},\\
    72(d^2-1), &\text{ if $d\ge5$ is odd}.
\end{cases}
\]
%For $d=2n\geq4$ one may take $|X_d|=12d(d+2)=48n(n+1)$, and for $d=2n+1\geq5$ one may take $|X_d|=72(d^2-1)=288n(n+1)$.
\end{theorem}
\begin{remark}
The absolute constants above can be improved, for instance, by using better constructions for Sidon sets, but we have chosen to keep them as is to keep the arguments simple.
The construction can be further improved when $n=p^k$ is a prime power because of the existence of mutually unbiased bases in dimension $n$. In this case one can construct $X_d$ with $|X_d|\le \tfrac32 d(d+2)=6n(n+1)$ in dimension $d=2n$ and $X_d$ with $|X_d|\le 9(d^2-1)$ in dimension $d=2n+1$.
\end{remark}

\section{Complex projective 2-designs with $O(n^2)$ points}

A finite nonempty multiset $Z \subset \Sb(\C^n)$, where $\Sb(\C^n)$ is the unit sphere in $\C^n$, is a \emph{complex projective $t$-design} if for every polynomial $p(z, \bar{z})$ that is homogeneous of degree $t$ in the coordinates $z_1, \ldots, z_n$ and homogeneous of degree $t$ in their complex conjugates $\bar{z}_1, \ldots, \bar{z}_n$, the following equality holds:
\[
    \frac{1}{|Z|} \sum_{z \in Z} p(z, \bar{z}) = \int_{\Sb(\C^n)} p(u, \bar{u}) \, d\mu(u),
\]
where $\mu$ is the normalized Haar measure on $\Sb(\C^n)$. Our key auxiliary result is the following lemma.

\begin{lemma}
    \label{lem:cp2design}
    For every $n\ge 2$ there exists a complex projective $2$-design of cardinality $8n(n+1)$.
\end{lemma}

For the proof, we require the concept of Sidon sets. Recall that a subset $A=\{a_1,\ldots,a_n\}\subseteq \Z/M\Z$ is called a Sidon set modulo $M$ if
	\[
        a_i+a_k\equiv a_j+a_l\pmod M
        \quad\Longleftrightarrow\quad
        \{i,k\}=\{j,l\} .
	\]
\begin{lemma}\label{lem:sidon}
	For every $n\geq2$ there exists a Sidon set of size $n$ modulo $M=8n^2$.
\end{lemma}
\begin{proof}
	Choose $q=2^l$ with $n\leq q<2n$.  The Bose--Chowla construction
	\cite{BoseChowla1962} gives a Sidon set $B$ of size $q$ modulo $q^2-1$. % $B\subset\Z/(q^2-1)\Z$ of size	$q$.
    Choose $n$ elements of $B$ and represent them by integers
	$a_1,\ldots,a_n\in\{0,\ldots,q^2-2\}$.  Since every sum $a_i+a_k$ lies in
	$[0,2q^2-4]$, and $2q^2-4<8n^2=M$, congruence modulo $M$ between two such sums
	implies equality as integers.  Equality as integers implies congruence modulo
	$q^2-1$, so the Sidon property modulo $q^2-1$ gives the required Sidon property modulo $M$. %in $\Z/(q^2-1)\Z$ gives the required Sidon property in $\Z/M\Z$.
\end{proof}

\begin{proof}[Proof of Lemma~\ref{lem:cp2design}.]
Let $A=\{a_1,\ldots,a_n\}\subset \Z/M\Z$ be the Sidon set from
Lemma~\ref{lem:sidon}.  Put $\omega=\exp(2\pi i/M)$ and, for
$r\in\Z/M\Z$,
\[
        f_r=\frac1{\sqrt n}\bigl(\omega^{ra_1},\ldots,\omega^{ra_n}\bigr)
        \in \C^n .
\]
Define the following labeled multiset of unit vectors:
\begin{equation}\label{eq:CPdesign}
        Z=
        \{f_r:r\in\Z/M\Z\}
        \sqcup
        \bigsqcup_{j=1}^n \{\underbrace{e_j,\ldots,e_j}_{8\text{ copies}}\}.
\end{equation}
Since $M=8n^2$, we have
\[
        |Z|=M+8n=8n(n+1).
\]
Next, we check the projective $2$-design condition.  It is enough to prove
that, for all $1\leq i,j,k,l\leq n$,
\begin{equation}\label{eq:cp2-moment}
        \frac1{|Z|}\sum_{z\in Z}
        z_i\overline {z_j} z_k\overline {z_l}
        =
        \frac{\delta_{ij}\delta_{kl}+\delta_{il}\delta_{jk}}{n(n+1)}.
\end{equation}
For the vectors $f_r$ we have
\[
\begin{aligned}
        \sum_{r\in\Z/M\Z}
        (f_r)_i\overline{(f_r)_j}(f_r)_k\overline{(f_r)_l}
        &=
        \frac1{n^2}
        \sum_{r\in\Z/M\Z}
        \omega^{r(a_i+a_k-a_j-a_l)}  \\
        &=
        \frac{M}{n^2}\,
        \bbone_{a_i+a_k\equiv a_j+a_l\pmod M}.
\end{aligned}
\]
By the Sidon property,
\[
        \bbone_{a_i+a_k\equiv a_j+a_l\pmod M}
        =
        \delta_{ij}\delta_{kl}+\delta_{il}\delta_{jk}
        -\delta_{i=j=k=l}.
\]
The eight copies of each coordinate vector give
\[
        \sum_{s=1}^n 8
        (e_s)_i\overline{(e_s)_j}(e_s)_k\overline{(e_s)_l}
        =8\,\delta_{i=j=k=l}.
\]
Using $M=8n^2$ and adding the two contributions, we obtain
\[
        \sum_{z\in Z}z_i\overline z_j z_k\overline z_l
        =8\bigl(\delta_{ij}\delta_{kl}+\delta_{il}\delta_{jk}\bigr).
\]
Dividing by $|Z|=8n(n+1)$ gives \eqref{eq:cp2-moment}.  Hence $Z$
defines a complex projective $2$-design in $\C\Pb^{n-1}$ (with repeating points).
\end{proof}

\section{Proof of Theorem~\ref{thm:main}}
For even dimensions $d=2n$ with $n\ge2$, we take the complex projective $2$-design $Z$ constructed in Lemma~\ref{lem:cp2design} and apply the fibration construction \cite[Thm 4.1]{Kuperberg2006}: under the identification of  $\C^n$ with $\R^{2n}$, $X_{2n}$ is obtained by taking $6$ points forming a regular hexagon on each fiber $\pi^{-1}(z)$ for $z\in Z$ where $\pi\colon\C^n\to\C\Pb^{n-1}$ is the canonical projection. Note that by choosing different rotations of the hexagon for repeating points of $Z$ we may arrange it so that $X_{2n}$ is a set. 
Then $X_{2n}$ is a $5$-design and
\[
        |X_{2n}|=6\cdot 8n(n+1)=48n(n+1).
\]
For $d=2$ take the regular hexagon on $\Sb^1$, which is a tight spherical $5$-design.

In odd dimensions we apply a standard lifting construction. Let $d\geq3$ be odd and define
	\[
        U_d=\Big\{\pm \alpha \colon \alpha^3-\frac{3}{d}\alpha^2+\frac{9}{d^2(d+2)}\alpha-\varepsilon=0\Big\},
	\]
where $\varepsilon>0$ is chosen to be sufficiently small so that $U_d$ is a subset of $[-1,1]$ (this can be done because for $\varepsilon\to 0$ the roots converge to $0$ and $\frac{3}{2d}(1\pm \sqrt{\frac{d-2}{d+2}})$).
Then the odd moments of $U_d$ vanish, and
	\[
        \frac1{6}\sum_{u\in U_d}u^2=\frac1d,
        \qquad
        \frac1{6}\sum_{u\in U_d}u^4=\frac3{d(d+2)}.
	\]
Thus $U_d$ has the moments of the projection of the normalized measure of $\Sb^{d-1}$ onto $[-1,1]$.

If $Y\subset \Sb^{d-2}\subset\R^{d-1}$ is a spherical $5$-design, set
\[
        \Lambda_d(Y)=
        \{(u,\sqrt{1-u^2}\,y):u\in U_d,
        \ y\in Y\}\subset \Sb^{d-1}.
\]
Then, by \cite[Thm. 4.1]{RabauBajnok1991}, $\Lambda_d(Y)$ is a spherical $5$-design.

For $d=2n+1\geq5$, take $Y=X_{2n}$ as constructed above, and set
$X_d=\Lambda_d(Y)$.  Then
\[
        |X_d|=6|X_{2n}|=288n(n+1)=72(d^2-1)<72d^2.
\]
For $d=3$, we may take $X_3$ to be the icosahedron, which is a tight spherical $5$-design of size $12$.  Together with $X_1=\{-1,1\}$, this proves Theorem~\ref{thm:main}.

\end{document}